\documentclass[oneside,reqno]{amsart}
\usepackage{amssymb}
\usepackage{amsmath}
\usepackage{amsthm}
\usepackage{amsbsy}
\usepackage{enumerate}
\usepackage{bm}
\usepackage{hyperref}
\date{\today}
\usepackage{cite}
\usepackage{array}
\usepackage{xcolor}
\usepackage{tikz}
\usepackage{graphics}
\usepackage{subcaption}
\usepackage{mathtools}
\usepackage{relsize}
\usepackage{soul}
\usepackage{mathdots}
\usepackage{multirow}
\usetikzlibrary{matrix}
\tikzset{highlight/.style={fill=red,draw=none,fill opacity=.2,inner sep=0pt}
}
\makeatletter
\newcommand{\distas}[1]{\mathbin{\overset{#1}{\kern\z@\sim}}}

\newcommand{\distras}[1]{
  \mathbin{\overset{#1}{\kern\z@\resizebox{\wd\mybox}{\ht\mysim}{$\sim$}}}
}

    \newtheorem{theorem}{Theorem}
    \newtheorem{lemma}[theorem]{Lemma}
    \newtheorem{proposition}[theorem]{Proposition}
    \newtheorem{corollary}[theorem]{Corollary}
    \newtheorem{claim}[theorem]{Claim}

\theoremstyle{definition} 
    \newtheorem{definition}[theorem]{Definition}

    \newtheorem{remark}[theorem]{Remark}
    \newtheorem{example}[theorem]{Example}
    \newtheorem{exercise}[theorem]{Exercise}
    
    \newtheorem{assumption}[]{Assumption}

\def\kap{\kappa}

\def\sig{\sigma}

\newcommand{\E}{\mathbb{E}}

\def\Z{\mathbb{Z}}

\def\N{\mathbb{N}}
\def\PP{{\mathcal P}}

\def\R{\mathbb{R}}

\def\WW{\mathcal{W}}

\def\tends{\rightarrow}

\def\<{\langle}
\def\>{\rangle}

\newcommand{\one}{{\mathbf 1}}

\newcommand{\Cov}{\operatorname{Cov}}

\newcommand\Tr{{\mbox{Tr}}}

\newcommand\mnote[1]{} 
\newcommand\be{\begin{equation*}}

\newcommand\ee{\end{equation*}}

\newcommand\ben{\begin{equation}}
\newcommand\een{\end{equation}}
\newcommand\bes{\begin{eqnarray*}}
\newcommand\ees{\end{eqnarray*}}

\newcommand\bex{\begin{exercise}}
\newcommand\eex{\end{exercise}}
\newcommand\beg{\begin{example}}
\newcommand\eeg{\end{example}}
\newcommand\benu{\begin{enumerate}}
\newcommand\eenu{\end{enumerate}}
\newcommand\beit{\begin{itemize}}
\newcommand\eeit{\end{itemize}}
\newcommand\berk{\begin{remark}}
\newcommand\eerk{\end{remark}}
\newcommand\bdefn{\begin{defintion}}
\newcommand\edefn{\end{definition}}
\newcommand\bthm{\begin{theorem}}
\newcommand\ethm{\end{theorem}}
\newcommand\bprf{\begin{proof}}
\newcommand\eprf{\end{proof}}
\newcommand\blem{\begin{lemma}}
\newcommand\elem{\end{lemma}}

\newcommand{\var}{\mbox{\rm Var}}
\newcommand{\sm}{{\raise0.3ex\hbox{$\scriptstyle \setminus$}}}

\def\sig{\sigma}

\def\tends{\rightarrow}

\def\CHI{\mathchoice%
{\raise2pt\hbox{$\chi$}}%
{\raise2pt\hbox{$\chi$}}%
{\raise1.3pt\hbox{$\scriptstyle\chi$}}%
{\raise0.8pt\hbox{$\scriptscriptstyle\chi$}}}
\def\smalloplus{\raise1pt\hbox{$\,\scriptstyle \oplus\;$}}

\usepackage{tikz}

\def\NN{{\mathcal N}}

\author[K. Kumar]{Kiran Kumar$^{\ast}$}
\address{$^{\ast}$Division of Science,
  New York University Abu Dhabi,
  Abu Dhabi, United Arab Emirates}
\email{kiran.kumar.as.math [at] gmail.com}

\author[S. N. Maurya]{Shambhu Nath Maurya$^{\dagger}$}
\address{$^{\dagger}$Department of Mathematics,
  Indian Institute of Technology Roorkee,
  247667, Uttarakhand, India}
\email{shambhumath4 [at] gmail.com}

\author[K. Saha]{Koushik Saha$^{\ddagger}$}
\address{$^{\ddagger}$Department of Mathematics, Indian Institute of Technology Bombay, Mumbai, 400076, India}
\email{koushik.saha [at] iitb.ac.in}

\date{\today}

\title[Time-dependent fluctuations for  generalized patterned random matrices] 
{A unified approach to fluctuations of LES of time-dependent generalized patterned random matrices}

 \usepackage{authblk}

\numberwithin{equation}{section}
\numberwithin{theorem}{section}
\begin{document}

\begin{abstract}
A unified approach for studying the linear eigenvalue statistics (LES) of generalized patterned random matrices was developed in \cite{KMS_UnifiedLES_24}. In this paper, as a continuation, we develop a unified method for studying the time-dependent LES of generalized patterned random matrices with independent process entries such as Brownian motion and Ornstein-Uhlenbeck process. We prove that if the limiting moment sequence of the LES of the generalized patterned random matrix with independent standard Gaussian entries is M-determinate, then under some minimal assumption on the link function, the time-dependent LES converges to a deterministic real-valued process.
\end{abstract}

\maketitle



\section{  Introduction} \label{sec:Intro}
  
Patterned  matrices such as Toeplitz, Hankel and circulant appear widely in different areas of mathematics. The study of their random versions was initiated in the review article by Z. D. Bai \cite{bai_zd_99}. In this work, we look at a generalization of patterned random matrices introduced in \cite{KMS_UnifiedLES_24}, which includes the banded and block versions of patterned random matrices.
For an input sequence $\{x_i : i \in \Z^d\}$ with $d \in \N$, an $N\times N$ generalized patterned random matrix is defined as
\begin{equation}\label{eqn:A_ij_general}
	(A)_{ij}=\begin{cases}
		x_{L(i, j)} &\text{ if } (i,j) \in \Delta, \\
		0 &\text{ otherwise},
	\end{cases}
\end{equation} 
where 
\begin{equation}\label{eqn: link_fun}
	L:\{1,2, \ldots, N\}^{2} \rightarrow \mathbb{Z}^{d}
\end{equation}
is a  symmetric link functions and $\Delta$ is a symmetric subset of $[N]^2$, i.e, $(i,j) \in \Delta$ implies $(j,i) \in \Delta$.
 Here $[N]$ denotes the set $\{1,2, \ldots,N\}$. Note that generalized patterned matrices are symmetric matrices.

 Here we consider a sequence of generalized patterned random matrices $\{A_n\}$, where the dimension of $A_n$ is $N(n)\times N(n)$ with a strictly increasing function $N:\mathbb N \to \mathbb N$. For notational convenience, we shall simply use $N$ in place of $N(n)$. 
 We have the following assumption on the link functions $\{L_n:[N(n)] \rightarrow \Z^d\}$ of $\{A_n\}$.

\begin{assumption}\label{Condition B}
	 There exists a constant $B < \infty$ such that
\begin{equation}\label{eqn:propertyB}
		\sup _n \sup _{t \in \mathbb{Z}^d} \sup _{1 \leq k \leq N} \#\{l: 1 \leq l \leq N, L_n(k, l)=t\} \leq B.
	\end{equation}
\end{assumption}

The study of \textit{linear eigenvalue statistics} (LES) is an important area of research in random matrix theory. For an $n \times n$ matrix $A$ and a test function $\phi$, the LES of $A$ is defined as
\begin{equation} \label{eqn:LES}
	\mathcal{A}(\phi)= \sum_{i=1}^{n} \phi(\lambda_i),
\end{equation} 
where $\lambda_1 , \lambda_2 , \ldots , \lambda_n$ are the eigenvalues of $A$. 

In \cite{KMS_UnifiedLES_24}, a unified approach was developed  for studying the LES of generalized patterned random matrices with independent input entries. The following theorem is a direct consequence of \cite[Theorem 2]{KMS_UnifiedLES_24}.

\begin{theorem}\label{thm:unified_moments_p_even}
	Let $\left\{A_n=\left(x_{L_n(i, j)} \one_{\Delta_n}(i,j)\right)_{i, j=1}^{N}\right\}$ be a sequence of symmetric generalized patterned random matrices where 
	\begin{enumerate}
		\item [(i)] the input sequence $\{x_i\}$ are independent random variables with mean zero, variance one, $\E(x_i^4)= c$ for all $i$ and uniformly bounded moments of all order,
		\item [(ii)] the sequence of link function $\{L_n\}$ satisfies Assumption \ref{Condition B}.
	\end{enumerate} 
	For even $p$, suppose $\displaystyle{\sigma_p^2 =\lim_{n \rightarrow \infty} \E\left(\frac{1}{\sqrt{N}} \left[\Tr\left(\frac{A_{n}}{\sqrt{N}}\right)^{p}-\mathbb{E} \Tr\left(\frac{A_{n}}{\sqrt{N}}\right)^{p}\right]\right)^2}$ exists, then
	\begin{align*}
		\frac{1}{\sqrt{N}} \left[\Tr\left(\frac{A_{n}}{\sqrt{N}}\right)^{p}-\mathbb{E} \Tr\left(\frac{A_{n}}{\sqrt{N}}\right)^{p}\right] \stackrel{d}{\longrightarrow} 
		\left\{\begin{array}{ll} 
			\displaystyle  N\left(0, \sigma_p^2\right)  & \text{if }       \sigma_p^2>0, \\
			\displaystyle 0 & \text{if } \sigma_p^2=0.	 
		\end{array}\right. 
	\end{align*}
\end{theorem}
Over the last decade, LES have been  investigated for various random matrix models with time-dependent input entries.  For Wigner matrices with Brownian motion entries, central limit theorem for LES was derived in \cite{LES_Wigner_brownian} and for Ornstein-Uhlenbeck process entries in \cite{LES_Wigner_OU}. CLT for the LES of Wigner matrices with more general Gaussian process entries was derived in \cite{LES_Wigner_Gaussian}. For band Toeplitz matrices with Brownian motion entries, fluctuations of LES were studied in \cite{li_sun_2015,m&s_toeplitz_2020}, and for circulant-type matrices and Hankel matrices with Brownian motion entries in \cite{rcsc_process} and \cite{kiran+maurya_hankel_22}.

In this article, we attempt to build a unified  framework to investigate the LES of time-dependent generalized patterned random matrices.
Let $\{a_i(t): t \geq 0, \ i \in \Z^d\}$ be a sequence of independent stochastic processes. We define the matrix-valued stochastic process $\{A_n(t); t\geq 0 , n\in \mathbb N\}$ corresponding to the link function $L_n$ and a subset $\Delta_n \subset [N]^2$, as  
$$ A_{n}(t)=\left(a_{L_n(i, j)}(t)\one_{\Delta_n}(i,j)\right)_{i,j=1}^N.$$
For the matrix-valued stochastic process $A_n(t)$,  
we define 
\begin{equation}\label{eqn:kappa_p(t)}
	\kappa_{p}(t) :=\frac{1}{\sqrt{N}} \left[\Tr\left(\frac{A_{n}(t)}{\sqrt{N}}\right)^{p}-\mathbb{E} \Tr\left(\frac{A_{n}(t)}{\sqrt{N}}\right)^{p}\right]. 
\end{equation}
Note that $\kap_p(t)$ is the time-dependent LES of $A_n(t)$ corresponding to the test function $\phi(x)=x^p$. For $t=1$, we denote $\kappa_p(t)$ by $\kappa_p$, i.e., 
\begin{equation}\label{eqn:w_cp}
	\kappa_{p} :=\frac{1}{\sqrt{N}} \left[\Tr\left(\frac{A_{n}(1)}{\sqrt{ N}}\right)^{p} -\mathbb{E} \Tr\left(\frac{A_{n}(1)}{\sqrt{ N}}\right)^{p}\right].
\end{equation}
Here we derive a unified result for generalized patterned random matrices with Brownian motion entries and Ornstein-Uhlenbeck process entries. More precisely, the objective is to derive process convergence for $\kappa_p(t)$ using asymptotic properties of the random variable $\kappa_{p}$.
 The main result (Theorem \ref{thm:process_brownian}) of this paper establishes that for both these processes,
\[
\big\{	\kappa_{p}(t) ; t \geq 0\big\}  \stackrel{\mathcal{D}}{\rightarrow} \big\{ \Gamma_{p}(t) ; t \geq 0\big\},
\]
provided that the limiting moment sequence of $\kappa_p$ is M-determinate. As an application of our main result, we prove several new results (Corollaries \ref{cor:process_LES_example_1} and \ref{cor:process_LES_example_2}) for banded and block patterned random matrices with Brownian motion entries and Ornstein-Uhlenbeck process entries. For patterned random matrices other than Wigner, this is the first work with Ornstein-Uhlenbeck process entries. Our work predominantly uses combinatorial techniques, in contrast to Ito calculus and martingale techniques usually employed in this context.
        
\vskip2pt
\noindent\textbf{Outline of paper:}
In Section \ref{sec:prelim}, we first discuss and state a few definitions and combinatorial sets which will be required to state our main result and later, to prove it. In Section \ref{sec:thm_Brown}, we prove our main result for Brownian  motion entries and in Section \ref{sec:thm_OH}, we establish the result for Ornstein-Uhlenbeck entries. 

\section{Preliminaries and main results} \label{sec:prelim}
\subsection{Random processes}
Consider the space $C_{\infty}:= C[0, \infty)$ of all real-valued continuous functions on $[0, \infty)$. It is well-known that $C_{\infty}$ is a metric space, and the Borel $\sig$-field generated by the open sets of $C_{\infty}$ is denoted by $\mathcal{C_{\infty}}$.
For a probability space $\left(\Omega,\mathcal{F},\mathbb{P}\right)$, we define a continuous process on $\Omega$ as a $C_\infty$-valued random variable. 

For real-valued continuous processes $\{\boldsymbol{X_n}\}_{n\geq 1}$ = $\{ X_n(t); t \geq 0 \}_{n\geq 1}$ and $\boldsymbol{X}=\{X(t);t\ge 0\}$, we say $\boldsymbol{X_n}$ converges to $\boldsymbol{X}$ in \textit{distribution}, if $  \E f(\boldsymbol{X_n}) \rightarrow  \E f(\boldsymbol{X}) \ \mbox{ as }n\to\infty,$
for every bounded, continuous real-valued function $f$ on $C_{\infty}$. We shall denote this convergence by $\boldsymbol{X_n} \stackrel{\mathcal D}{\rightarrow} \boldsymbol{X}$. 

In this paper, we consider the following two processes:

\begin{definition}\label{defn:BM}
A stochastic process $\{W_t\}_{t\geq 0}$ is called a standard Brownian
motion if:
\begin{enumerate}
    \item[(i)] $W_0 = 0$ almost surely.
    \item[(ii)] It has independent increments.
    \item[(iii)] For every $0 \leq s < t$,
    $
        W_t-W_s \stackrel{D}{\sim} \mathcal{N}(0,t-s).
    $
    \item[(iv)] Its sample paths $t \mapsto W_t$ are continuous almost surely.
\end{enumerate} 
\end{definition}

\begin{definition}\label{defn:OU}
  A stochastic process $\{X_t\}_{t\geq 0}$ is called a centered Ornstein--Uhlenbeck process if
  \begin{enumerate}
      \item[(i)] $X_t$ is stationary Gaussian with $X_t \stackrel{D}{\sim} \NN(0,\frac{\sigma^2}{2\theta})$.
      \item[(ii)] $X_t$ is Markovian, i.e. for all $t_1 \leq t_2 \leq \cdots \leq t_k$ and $x \in \R$,
      \[
      \mathbb{P}[X_{t_k} \leq x| X_{t_1},X_{t_2},\ldots, X_{t_{k-1}}]= \mathbb{P}[X_{t_k} \leq x|  X_{t_{k-1}}].
      \]
      \item[(iii)] Its sample paths $t \mapsto X_t$ are continuous almost surely.
  \end{enumerate}  
\end{definition}
\subsection{Sentences, clusters and circuits} \label{sec:words}
In this section, we recall some basic combinatorial objects from \cite[Section 2]{KMS_UnifiedLES_24} required to state our main results. 
\begin{definition}\label{defn:word}
	Let $[p]=\{1,2,\ldots,p\}$. We can label partitions of $[p]$ by words of length $p$ such that the first occurrence of each letter is in alphabetical order. 
    For example, for $p=5$, the partition $\{\{1,3,5\},\{2,4\}\}$ is represented by the word \textit{ababa}. 
    
    The set of all words of length $p$ is denoted by $\WW_p$. A word $w \in \WW_p$ is called a pair-partition if each letter in $w$ repeats exactly twice. The set of all pair-partitions of $[p]$, is denoted by $\mathcal{P}_2(p)$.
\end{definition}

\begin{definition}\label{def:sentence}
	For $k \geq 1$, a partition of $\{1,2,\ldots, k\} \times \{1,2,\ldots,p\}$ is represented by a sentence with $k$ words each of length $p$, by assigning different letters to different blocks of the partition so that the first occurrence of the letters in the sentence is in the alphabetical order. In particular, the $i$-th word of the sentences represent the elements in $\{i\} \times [p]$. We denote the set of sentences with $k$ words each of length $p$ is by $\WW_{p,k}$.
\end{definition}
\begin{definition}
For a word $w \in \WW_p$, we define by $S_w$ the multi-set of letters in $w$. For a sentence $W=(w_1,w_2,\ldots ,w_k) \in \WW_{p,k}$, we define $S_W$ as the multi-set of all letters in $W$.
\end{definition}

\begin{definition}\label{defn: P_2,2 , P2,4}
	We consider the following type of sentences:
	
		\noindent (a) For a positive integer $p$, we denote a subset of $\WW_{p,2}$ by $\mathcal{P}_2(p, p)$, which consists of sentences $W=(w_1,w_2) \in \WW_{p,2}$ such that each letter in $W$ appears exactly twice and $S_{w_1} \cap S_{w_2} \neq \phi$.
		\vskip2pt
		\noindent (b) For even positive integer $p$, we denote by $\mathcal{P}_{2,4}(p, p)$, the set of all sentences $W=(w_1,w_2) \in \WW_{p,2}$ for which there exists $x \in S_W$ such that
		\begin{enumerate}
        \item [(i)] $x$ appears exactly four times in $S_W$,
		\item [(ii)]the multi-set intersection of $S_{w_1}$ and $S_{w_2}$ is $\{x,x\}$,
        \item [(iii)] all other letters appear exactly twice in $S_W$.
		\end{enumerate}
		For odd $p$, we define $\mathcal{P}_{2,4}(p, p)$ as an empty set.
\end{definition} 

Now, we define the notion of a graph associated with sentences.

\begin{definition}\label{defn_cluster}
	Let $W=(w_1,w_2,\ldots,w_k) \in \mathcal{W}_{p,k}$ be a sentence. We  construct a graph $G_{W}=(V_{W},E_{W})$ associated with $W$ in the following fashion: We define $V_{W}=\{ 1, 2,\ldots , k\}$ and for $r \neq s$, $\{r,s\} \in E_{W}$ if $S_{w_r} \cap S_{w_s} \neq \phi$. The vertex set of a maximal connected subgraph of $G_W$ is called a \textit{cluster}.
\end{definition}

\begin{definition}
   For $p \geq 1$ and $k \geq 3$, $W \in \mathcal{W}_{p,k}$ is called a \textit{clique sentence} if there exists a letter $x \in S_W$ such that
	\begin{enumerate}
    \item [(i)] $S_{w_i} \cap S_{w_j}=\{x\}$ for all $i \neq j$, and $x$ appears exactly $k$ times in $S_W$,
	\item[(ii)] all other letters appear exactly twice in $S_W$.
	\end{enumerate} 
\end{definition}
\begin{definition}\label{def: special partition}
	For $p \geq 1$ and $k \geq 3$, a sentence $W \in \mathcal{W}_{p,k}$ is called a \textit{special partition} if each cluster $C$ of $G_W$ is of size at least two and
	\begin{enumerate}
		\item[(i)] if $|C| =2$, then $W_C$ belongs to either $\PP_2(p,p)$ or $\PP_{2,4}(p,p)$, and
		\item [(ii)] if $|C| \geq 3$, then $W_C$ is a clique sentence,
	\end{enumerate}
	where $W_C:=(w_i)_{i \in C}$.
	The set of all special partitions in $\mathcal{W}_{p,k}$ is denoted by $SP(p,k)$.
\end{definition}
Next, we define the notion of circuits.
Recall the definition of link function $L: [N] \rightarrow \Z^d$ from \eqref{eqn: link_fun}. 
\begin{definition}\label{defn:circuit}
	Consider a link function $L: [N] \rightarrow \Z^d$, where $d \in \N$. A function $\pi:\{0,1,2, \ldots, p\} \rightarrow[N]$ with $\pi(0)=\pi(p)$ is called a \textit{circuit} of length $p$. Two circuits of length $p$, $\pi_1$ and $\pi_2$ are said to be equivalent if their $L$-values match at the same locations. That is, for all $1 \leq i, j \leq p$, 
	\begin{align*}
		 L(\pi_1(i-1), \pi_1(i)) =L(\pi_1(j-1), \pi_1(j)) 
		\Longleftrightarrow  L(\pi_2(i-1), \pi_2(i)) =L(\pi_2(j-1), \pi_2(j)).
	\end{align*}
\end{definition}
The above condition gives an equivalence relation on the set of all circuits of length $p$. It is immediate that any equivalence class can be indexed by a partition of $\{1,2, \ldots, p\}$ where each block of the given partition identifies the positions where the $L$-matching take place. We can label these partitions by words of length $p$ such that the first occurrence of each letter is in alphabetical order. For example, if $p=5$, then the partition $\{\{1,3,5\},\{2,4\}\}$ is represented by the word \textit{ababa}. This partition identifies all circuits $\pi$ for which $L(\pi(0), \pi(1))=L(\pi(2), \pi(3))=L(\pi(4), \pi(5))$ and $L(\pi(1), \pi(2))=L(\pi(3), \pi(4))$.

For a link function $L$ and a sentence $W=(w_1,w_2,\ldots ,w_k ) \in \WW_{p,k}$, we define 
\begin{equation*}
	\Pi_L(W)=\{(\pi_1,\pi_2,\ldots ,\pi_k ): w_r[i]=w_s[j] \Leftrightarrow L(\pi_r(i-1), \pi_r(i))=L(\pi_s(j-1), \pi_s(j))\},
\end{equation*}
where $w_r[i]$ denotes the $i$-th letter of $w_r$. Similarly, for a set $\Delta \subset [N]^2$, we define the set 
\begin{align*}
	\Pi_{L,\Delta}(W) &=\{(\pi_1,\pi_2,\ldots ,\pi_k ) \in \Pi_{L}(W) :  (\pi_u(i-1), \pi_u(i)) \in \Delta \ \forall \ 1 \leq u \leq k , 1 \leq i \leq p \}. 
\end{align*}
Note that in the above equation, the circuits are maps from $\{0,1,\ldots, p\}$ to $[N]$ and hence the cardinality of $\Pi_{L,\Delta}$ depends on $N$ as well. To simplify the notations, we are suppressing this dependence in the notation.

\subsection{Main Results}
Now we state our main result. Recall $\kappa_p(t)$ and $\kappa_p$ defined in \eqref{eqn:kappa_p(t)} and \eqref{eqn:w_cp}, and $SP(p,k)$ from Definition \ref{def: special partition}.
\begin{theorem} \label{thm:process_brownian}
Let $p \in \N$ and $\{A_{n}(t)=\left(a_{L_n(i, j)}(t)\one_{\Delta_n}(i,j)\right)_{i,j=1}^N\}$ be a sequence of matrix-valued stochastic processes such that 
\begin{enumerate}
	\item [(i)]$\{a_i(t)\}$ are either independent standard Brownian motions or independent centered Ornstein-Uhlenbeck processes. 
	\item [(ii)] the sequence of link function $\{L_n\}$ satisfies Assumption \ref{Condition B},
	\item [(iii)]for each sentence $W \in SP(p,k)$, the limit $\theta(W)=\lim_{n \rightarrow \infty} \frac{1}{N^{\frac{pk+k}{2}}}\#\Pi_{L_n,\Delta_n}\left(W\right) $ exists,
	\item [(iv)]for a positive integer $p$, the limiting moment sequence of $\kappa_p$ is \textit{M-determinate}.
\end{enumerate} 
Then there exists a real-valued continuous process $ \big\{ \Gamma_{p}(t) ; t \geq 0\big\}$ such that
	\begin{equation}
		\big\{	\kappa_{p}(t) ; t \geq 0\big\}  \stackrel{\mathcal{D}}{\rightarrow} \big\{ \Gamma_{p}(t) ; t \geq 0\big\}.
	\end{equation}
  
Moreover, if the limiting distribution of $\kappa_p$ is Gaussian, then the limiting process  $ \big\{ \Gamma_{p}(t) ; t \geq 0\big\}$ is a Gaussian process.
\end{theorem}

Using Theorem \ref{thm:process_brownian}, we prove a host of new results for different classes of patterned random matrices.
In the first corollary, we consider the classical patterned random matrices (CPRMs), namely, Toeplitz, Hankel, symmetric circulant, and reverse circulant matrices, and their banded versions. We say that a matrix $A$ is banded with bandwidth $b_n$ if $A_{ij}=0$ for all $|i-j| \geq b_n$ and anti-banded if $A_{ij}=0$ for all $|i+j-(n+1)| \geq b_n$.
\begin{corollary}\label{cor:process_LES_example_1}
	Let $\{C_n(t)\}$ be a sequence of banded or anti-banded CPRMs with independent Brownian motion entries or independent Ornstein-Uhlenbeck processes with bandwidth $b_n$ such that $b_n/n \to b>0$. Then the time-dependent LES of $C_n(t)$ converges to a Gaussian process for all even degree monomial test functions.
\end{corollary}
Recall that block patterned random matrices fall into the class of generalized patterned random matrices. We consider block patterned matrices of the form:
\[A=(B_{L_1(i,j)})_{i,j=1}^M \text{ with } B_\ell= (x_{\ell,L_2(i,j)})_{i,j=1}^m \text{ for } \ell \in \operatorname{Range}(L_1),\]
where $L_1, L_2$ are appropriate link functions (for details, see \cite[Section 5]{KMS_UnifiedLES_24}). Note that the actual dimension of $A$ is $Mm \times Mm$, where $M^2$ is the number of blocks and each block is of size $m \times m$.

\begin{corollary}\label{cor:process_LES_example_2}
	Let $\{A_n(t)\}$ be a sequence of block matrices of one of the following classes with independent Brownian motion entries or independent Ornstein-Uhlenbeck processes:
	\begin{enumerate}
		\item [(a)] Block Wigner matrices with $M$ fixed, and blocks are CPRMs.
		\item [(b)] Block versions of CPRMs with Wigner blocks of fixed size $m\times m$.
		\item [(c)] Block versions of CPRMs, whose blocks are any one of the CPRMs.
	\end{enumerate}
	For all of the above matrices, the time-dependent LES converges to a Gaussian process for all even degree monomial test functions.
\end{corollary}
We skip the detailed proof of the above corollaries. Note that for classical patterned random matrices and their band versions, condition (iii) of Theorem \ref{thm:process_brownian} follows from \cite[Theorem 6.1]{KMS_UnifiedLES_24}, and for block patterned random matrices condition (iii) follows as a consequence of \cite[Theorem 6.1]{KMS_UnifiedLES_24} and \cite[Proposition 5.1]{KMS_UnifiedLES_24}. Lastly, condition (iv) is satisfied by Theorem \ref{thm:unified_moments_p_even}.
\section{Proof of Theorem \ref{thm:process_brownian}: Brownian motion entries} \label{sec:thm_Brown}
In Section \ref{subsec:proof_idea}, we outline the general framework underlying the proof of convergence of the time-dependent LES. In Sections \ref{subsec:finite_dim} and \ref{subsec:tight}, we establish our main result in two steps, proving, respectively, finite-dimensional convergence and tightness. Finally, in Section \ref{subsec:gaussian}, we consider the special case in which the limiting process is Gaussian.
\subsection{Proof Idea}\label{subsec:proof_idea}

Suppose $\{\boldsymbol{X_n}\}_{n\geq 1}$ = $\{ X_n(t); t \geq 0 \}_{n\geq 1}$ is a sequence of continuous processes on a probability space $(\Omega,\mathcal F,\mathbb P)$.
The following result follows from \cite[Theorem I.4.2]{Ikeda1981}.
\begin{lemma}
\label{lem:two conditions} 
	Suppose the sequence of continuous processes $\{\boldsymbol{X_n}\}$ on $(\Omega,\mathcal F,\mathbb P)$ obeys the following two conditions
	\begin{enumerate}
		\item[(i)]  \textbf{Finite dimensional convergence:} For all $r \in \mathbb{N}$ and $0\leq t_1<t_2< \cdots <t_r$, as $n \tends \infty$,
		$(X_{n}(t_1), X_{n}(t_2), \ldots , X_{n}(t_r))$ converges in distribution.
		\vskip2pt
		\item[(ii)] 	 \textbf{Tightness:} The sequence  $\{\boldsymbol{X_n}\}_{n\geq 1}$ is tight. 
	\end{enumerate}
	Then there exists a probability space $(\widehat{\Omega}, \widehat{\mathcal{F}},\widehat{\mathbb{P}})$ and continuous processes $\boldsymbol{\widehat{X}_n}, \widehat{\boldsymbol{X}}$ on $\widehat{\Omega}$ such that $\boldsymbol{\widehat{X}_n}  \stackrel{D}{\sim}\boldsymbol{X_n}$ for each $n$ and $\boldsymbol{\widehat{X}_n}$ converges to $\widehat{\boldsymbol{X}}$ almost surely, where $Y \distas{D} Z$ means that $Y$ and $Z$ have the same distribution. 
\end{lemma}

Since $\widehat{\boldsymbol{X}}$ is a continuous process, it defines a probability measure $\mu_{\widehat{X}}$ on $( C_{\infty},\mathcal{C_{\infty}} )$. Here, we say that $\boldsymbol{X_n}$ converges in distribution to $\mu_{\widehat{X}}$. The following result provides the existence of a continuous process on $\Omega$ with distribution $\mu_{\widehat{X}}$.

\begin{lemma} \cite[Lemma 2]{Berti_Luca_Rigo}\label{res:exist of process}
	Suppose $(\Omega, \mathcal{F}, \mathbb{P})$ is a probability space,  $\mathcal{D}$ is a Polish space and $X_n: \Omega \rightarrow \mathcal{D}$ converges in distribution to a probability measure $\mu$ on $\mathcal{D}$. Then 
    there is a measurable map $X : \Omega \rightarrow \mathcal{D}$ with distribution $\mu$. 
\end{lemma}

It is well-known that $C_{\infty}$ is a \textit{Polish space}. Therefore, it follows from Lemma \ref{res:exist of process} that for any sequence $\{ \boldsymbol{X_n}\}$ on $\Omega$ obeying the conditions (i) and (ii) of Proposition \ref{lem:two conditions} there exists a continuous process $\boldsymbol{X}$ on $\Omega$ such that $\boldsymbol{X_n}\stackrel{D}{\rightarrow} \boldsymbol{X}$. 

To prove Lemma \ref{lem:two conditions} for the process $\kap_p(t)$, we prove the following two propositions. Proposition \ref{pro:finite convergence} gives the finite-dimensional convergence of $\kap_p(t)$ and, as a consequence of \cite[Theorem I.4.3]{Ikeda1981}, Proposition \ref{pro:tight} gives the tightness of $\kappa_p(t)$.

\begin{proposition} \label{pro:finite convergence}
	Let $p$ be a positive integer. Then for any $ r \in \mathbb{N}$ and $0\leq t_1<t_2< \cdots <t_r$, as $n \tends \infty$,
	$(\kappa_{p}(t_1), \kappa_{p}(t_2), \ldots , \kappa_{p}(t_r))$ converges in distribution.
\end{proposition}

\begin{proposition} \label{pro:tight}
	For each positive integer $p$, there exist positive constants $C$ and $\gamma$ such that
	\begin{equation*} \label{eq:tight3}
		\E{|\kap_{p}(0)|^ \gamma } \leq C \ \ \ \forall \ n \in \mathbb{N}. 	 
	\end{equation*}
	Furthermore, there exist positive constants $\alpha, \beta $ and $M_T$, $T= 1, 2, \ldots$ such that 
	\begin{equation} \label{eqn:tight4}
		\E{|\kap_{p}(t)- \kap_{p}(s)|^ \alpha } \leq M_T |t-s|^ {1+ \beta} \ \ \ \forall \ n \in \mathbb{N} \ \mbox{and } t,s \in[0,\ T],(T= 1, 2, \ldots). 
	\end{equation}
\end{proposition} 

In Section \ref{subsec:finite_dim}, we prove Proposition \ref{pro:finite convergence}, and in Section \ref{subsec:tight} we prove Proposition \ref{pro:tight}.

\subsection{Finite Dimensional Convergence} \label{subsec:finite_dim}
 
Our approach in proving Proposition \ref{pro:finite convergence} is through the \textit{multi-dimensional moment method} (Lemma \ref{res:mutimoment_method}).  In Lemma \ref{pro:lim+E[eta]finite}, we show that $\E [q(\kappa_{p}(t_1),\ldots , \kappa_{p}(t_r))]$ converges for all monomials $q \in \R[x_1,\ldots , x_r]$ in $r$ variables. Then using Haviland's theorem (Lemma \ref{res:haviland_existance}), we establish that for all choices of $0 \leq t_1 <t_2 < \cdots <t_r$, the limits of $\E [q(\kappa_{p}(t_1),\ldots , \kappa_{p}(t_r))]$ determine a unique distribution $\Gamma_{t_1,\ldots t_r}$ on $\R^r$ (Lemma \ref{lem:multi_dim}). Finite-dimensional convergence of $\kappa_p(t)$ then follows from Lemma \ref{res:mutimoment_method}.
\begin{lemma}  \cite[Theorem 1]{haviland_multimoment}\label{res:mutimoment_method}
	Let $\mu_n$ and $\mu$ be probability measures on $\mathbb{R}^r$. Suppose
	that $\mu$ is uniquely determined by its moments and
	\begin{equation} \label{eqn:multimoment}
		\lim_{n\to\infty}	\int_{\mathbb{R}^r} x_1^{k_1}x_2^{k_2}\cdots x_{r}^{k_r} \mu_n(dx_1, \ldots, dx_r) = 	\int_{\mathbb{R}^r} x_1^{k_1}x_2^{k_2}\cdots x_{r}^{k_r} \mu(dx_1, \ldots, dx_r),
	\end{equation}
	for all $k_1, k_2, \ldots, k_r \in \{0,1,2, \ldots\}$. Then $\mu_n \stackrel{d}{\rightarrow} \mu $.
\end{lemma}

For a fixed $r \in \N$ and $0 \leq t_1 < t_2< \cdots < t_r$, consider the $r$-dimensional random vector $(\kap_{p}(t_1), \kap_{p}(t_2), \ldots , \kap_{p}(t_r))$.
In the following lemma we show that $\E\left[q\left(\kap_p(t_1),\ldots,\kap_p(t_r)\right)\right]$ is convergent for every  monomial $q$.
\begin{lemma} \label{pro:lim+E[eta]finite}
	Let $A_n(t)=(a_{L_n(i,j)}(t) \one_{\Delta_n}(i,j))_{i,j=1}^N$ be a symmetric matrix-valued stochastic process, where $\{a_i(t)\}$ are independent Brownian motions, obeying conditions (ii)-(iii) of Theorem \ref{thm:process_brownian}.
  Then for every choice of  $0 \leq t_1 < t_2 < \cdots < t_r$ and monomial $q \in \R [x_1,x_2,\ldots, x_r ]$, $\E\left[q\left(\kap_p(t_1),\ldots,\kap_p(t_r)\right)\right]$ converges.
\end{lemma}
\begin{proof} 
For $k_1,k_2,\ldots,k_r \in \N$, consider an arbitrary monomial $q(x_1,x_2,\ldots
,x_r)=x_1^{k_1}x_2^{k_2}\cdots x_r^{k_r}$.  Recall the definition of $\kap_p(t)$ from (\ref{eqn:kappa_p(t)}) and by expanding the trace, we get that
\begin{multline}\label{eqn: kappa p power k}
    	\kap_{p}^{k_i}(t_i)= \frac{1}{N^{\frac{pk_i+k_i}{2}}}\displaystyle \sum_{\pi_1,\ldots , \pi_k} \big[(a_{\pi_1}(t_i)- \E a_{\pi_1}(t_i))(a_{\pi_2}(t_i)- \E a_{\pi_2}(t_i)) \cdots (a_{\pi_k}(t_i)- \E a_{\pi_k}(t_i)) \\\prod_{u=1,2,\ldots , k \atop j=1,2,\ldots , p} \one_{\Delta_n}(\pi_u(j-1),\pi_u(j))\big],
\end{multline}
where for $u=1,2, \ldots,k_i$; $\pi_u:\{0,1,2,\ldots, p\} \rightarrow \{1,2,\ldots, N\}$ are circuits of length $p$ and $a_{\pi_u}$ is defined as  
\[
a_{\pi_u}=a_{L_n(\pi_u(0),\pi_u(1))}a_{L_n(\pi_u(1),\pi_u(2))}\cdots a_{L_n(\pi_u(p-1),\pi_u(p))}.\]
Since $\{a_i\}$ are i.i.d., grouping together all the terms in the equivalence class as per Definition \ref{def:sentence}, we have that
	\begin{equation}\label{eqn: E eta_p time product}
		\mathbb{E}[	\kap_p^{k_1}\left(t_1\right) \cdots \kap_p^{k_r}\left(t_r\right)] 
		= \frac{1}{N^{ \frac{pk+k}{2}}}  \sum_{W \in \WW_{p,k}} \sum_{\Pi_{L_n,\Delta_n}(W)}  \mathbb{E}\left(\prod_{i=1}^k \big(a_{\pi_i} (s_i)-\mathbb{E} a_{\pi_i}(s_i) \big)\right),
	\end{equation}
	where $k=k_1+k_2+\cdots +k_r$ and the summation is over $(\pi_1,\pi_2,\ldots, \pi_k) \in \Pi_{L_n, \Delta_n}(W)$, and for any $m \in \{1,2, \ldots,r\}$, 
    \begin{equation}
        s_i=t_m \quad \text{ for } \quad k_1+k_2+\cdots +k_{m-1}<i \leq k_1+k_2+\cdots+ k_m.\label{eqn:si}
    \end{equation}
	Note that since $a_i(t)$ are i.i.d., $\mathbb{E}\left(\prod_{i=1}^k \left(a_{\pi_i} (s_i)-\mathbb{E} a_{\pi_i}(s_i)\right)\right)$ in (\ref{eqn: E eta_p time product}) depends only on the sentence $W=(w_1,w_2,\ldots,w_k)$ and is fixed for every choice of $(\pi_1, \pi_2, \ldots , \pi_k) \in \Pi_{L_n,\Delta_n}(W)$ and we denote it by $\mathbb{E} \left(\prod_{i=1}^k \left(a_{w_i} (s_i)-\mathbb{E}a_{w_i}(s_i)\right)\right)$.  Further, we denote 
	\begin{equation}\label{eqn: v(W) defn}
		v(W)= \frac{\# \Pi_{L_n, \Delta_n}(W)}{N^{\frac{pk+k}{2}}}\mathbb{E}\left(\prod_{i=1}^k (a_{w_i} (s_i)-\mathbb{E} a_{w_i}(s_i))\right).
	\end{equation}
    Note that $\mathbb{E}\left(\prod_{i=1}^k (a_{w_i} (s_i)-\mathbb{E} a_{w_i}(s_i))\right)$ is a constant independent of $n$. Therefore,  it follows from \cite[Proposition 20]{KMS_UnifiedLES_24} that $v(W)$ is of the order $o(1)$ if $W \notin SP(p,k)$. Further, if $W \in SP(p,k)$, by condition (iii) of Theorem \ref{thm:process_brownian}, $\lim_{n \to \infty} v(W)$ exists.

For later use, we compute $\mathbb{E}\left(\prod_{i=1}^k (a_{w_i} (s_i)-\mathbb{E} a_{w_i}(s_i))\right)$, for sentences $W \in SP(p,k)$.  We have that for all clusters $C$ of $W$, $W_C$ is an element of $\PP_{2}(p,p)$ or $\PP_{2,4}(p,p)$ or a clique sentence. Direct computation using Lemma 4.5 of \cite{Marcus_ray_brownian_moments_1992} yields that 	 
\begin{align}
	& \mathbb{E}\left[\prod_{i  \in C} (a_{w_i} (s_i)-\mathbb{E} a_{w_i}(s_i))\right] \label{eqn:contri_clusters} \nonumber	 \\
	 &	=\begin{cases}
	 		\displaystyle{\sum_{\pi \in \PP_2(\ell)}\prod_{\{u,v\} \in \pi}\min \{s_{i_u},s_{i_v}\}} \prod_{u=1}^\ell(s_{i_u})^{\frac{p-1}{2}} & \text{ if } W_C=(w_{i_1},\ldots, w_{i_\ell}) \text{ is  a clique sentence},\\
	 	s_{i_1}^{\frac{p+c_0}{2}} s_{i_2}^{\frac{p-c_0}{2}} &	\text{ if } W_C=(w_{i_1},w_{i_2}) \in \PP_{2}(p,p), \\
	 		s_{i_1}^{\frac{p}{2}-1} 
        s_{i_2}^{\frac{p}{2}-1}(2s_{i_1}+s_{i_2}-s_{i_1}s_{i_2}) &	\text{ if } W_C=(w_{i_1},w_{i_2}) \in \PP_{2,4}(p,p),
	 \end{cases}
\end{align}
where for $(w_{i_1},w_{i_2}) \in \PP_{2}(p,p)$, $c_0$ is the number of  elements in $S_{w_{i_1}} \cap S_{w_{i_2}}$.
\end{proof}
Now using Haviland's theorem we show that the limits of sequence of joint moments determine a distribution. The uniqueness of this distribution is established by using its connection to the one-dimensional determinate problem.

\begin{lemma} \cite[Theorem 3.1.2]{Murray_Marshall} (\textit{Haviland's Theorem}) \label{res:haviland_existance}
For a linear function \\ $L: \R [x_1,x_2,\ldots, x_r ] \rightarrow \R$, and a closed set $K$ in $\R^r$, the following are equivalent:
	\begin{enumerate}
		\item [(i)]there exists a probability measure $\mu$ on $K$ such that 
		\begin{equation}
			L(f)=\int f\  \mathrm{d}\mu \quad \forall \ f \in \R [x_1,x_2,\ldots, x_r ].
		\end{equation}
		\item [(ii)]$L(f) \geq 0$, for all $f \in \R [x_1,x_2,\ldots, x_r ]$ such that $f \geq 0$ in $K$ and $L(\mathbf{1})=1, $ where $\mathbf{1}$ is the constant function with range $\{1\}$.
	\end{enumerate}
\end{lemma}
The following theorem provides the existence and uniqueness of a distribution whose moments are the limiting sequence of joint moments of the process $\{\kap_p(t): t \geq 0 \}$.
\begin{lemma}\label{lem:multi_dim}
Let $A_n(t)=(a_{L_n(i,j)}(t) \one_{\Delta_n}(i,j))_{i,j=1}^N$ be a symmetric matrix-valued stochastic process, where $\{a_i(t)\}$ are independent Brownian motions, obeying conditions (ii)-(iv) of Theorem \ref{thm:process_brownian}.
 Then for all $0 \leq t_1 < t_2 < \cdots < t_r$, there exists a unique distribution $\Gamma_{t_1,t_2,\ldots , t_r}$ in $\R^r$ such that for all $q \in \R [x_1,x_2,\ldots, x_r ]$ 
	\begin{equation}
		\int_{\R^r} q(x_1,\ldots, x_r) \,\mathrm{d}\Gamma_{t_1,t_2,\ldots , t_r}= \lim_{n \rightarrow \infty} \E\left[q\left(\kap_p(t_1),\ldots,\kap_p(t_r)\right)\right].
	\end{equation}
\end{lemma}
\begin{proof}
	We prove the existence part using Lemma \ref{res:haviland_existance}. Recall that $\kap_p(t_r)$ depends on $n$. For a fixed $n$, let $F_n$ denote the distribution of $\left(\kap_p(t_1),\kap_p(t_2), \ldots, \kap_p(t_r)\right)$ in $\R^r$. 
	
Consider the sequence of linear maps $L_n: \R[x_1,x_2,\ldots, x_r] \rightarrow \R$ defined by
	$$L_n(q):= \int q(x_1,x_2,\ldots, x_r)\mathrm{d}F_n.$$
By Haviland's theorem, we obtain that $L_n(q) \geq 0$ for all $q \in \R[x_1,x_2,\ldots , x_r]$ such that $q \geq 0$ on $\R^r$.

Now for $q \in \R[x_1,x_2,\ldots,x_r]$, we define $L(q) := \lim_{n \rightarrow \infty} L_n(q).$
Clearly, $L$ is a linear map, and $L(q) \geq 0$ for all $q \in \R[x_1,x_2,\ldots , x_r]$ such that $q \geq 0$ on $\R^r$. Further, since $F_n$ are probability distributions, we have $L_n(\mathbf{1})=1$ and therefore we obtain that $L(\mathbf{1})=1$. Hence by Lemma \ref{res:haviland_existance}, we get the existence of a probability measure.

Note from Theorem 3 of \cite{Petersen_carleman_multi} that to prove the uniqueness of $\Gamma_{t_1,t_2,\ldots , t_r}$, it is sufficient to show that for each $i=1,2,\ldots ,r$, the moment sequence $\{m_k^{(t_i)}= \E[\kap_p(t_i)]^k\}$ uniquely determines a measure on $\R$. Note that since $\{a_i(t)\}$ are i.i.d. Brownian motion entries, $\frac{a_i(t)}{\sqrt{t}}$ has the standard normal distribution. This shows that $\frac{\kap_p(t_i)}{(\sqrt{t_i})^p} \stackrel{D}{\sim} \kappa_p$, where $\kappa_p$ is defined in (\ref{eqn:w_cp}). As the limiting moment sequence of $\kap_p$ is M-determinate,  $\{\frac{m_k^{(t_i)}}{t_i^{p/2}}\}$ defines a unique measure on $\R$, which in turn implies that $\{m_k^{(t_i)}\}_{k \in \N}$ determines a unique measure on $\R$.
\end{proof} 
It follows from Lemma \ref{lem:multi_dim} and Lemma \ref{res:mutimoment_method} that $\left(\kap_p(t_1),\ldots,\kap_p(t_r)\right)$ converges in distribution to $\Gamma_{t_1,t_2,\ldots , t_r}$ for every choice of $ 0 \leq t_1 < t_2 < \cdots < t_r$.
\subsection{Tightness} \label{subsec:tight}
\begin{proof}[Proof of Proposition \ref{pro:tight}] First note from (\ref{eqn:kappa_p(t)}) that $\kap_p(0)=0$, and  to complete the proof, we prove (\ref{eqn:tight4}) for $\alpha = 4$ and $\beta = 1$. 
	Suppose $p$ is a positive integer and $t,s \in[0,T] $, for some fixed $T \in \mathbb{N}$. Then from (\ref{eqn:kappa_p(t)}),
	\begin{align} \label{eqn:w_p(t-s)}
		\kap_{p}(t) - \kap_{p}(s) &= \frac{1}{N^{\frac{p+1}{2}}} \Big[ \Tr(A_n(t))^{p} - \Tr(A_n(s))^{p} - \E[\Tr(A_n(t))^{p} - \Tr(A_n(s))^{p}]  \Big].
	\end{align}
	For $0< s <t $, using binomial expansion, we get    
	\begin{equation} \label{eqn:trace t-s}
		\Tr (A_n(t))^{p} - \Tr (A_n(s))^{p}  = \Tr[ (A'_n(t,s) )^{p}] + \sum_{d=1}^{p-1} \binom{p}{d} \Tr [(A'_n(t,s))^d (A_n(s))^{p-d}],
	\end{equation}
	where $A'_n(t,s)=\left((a_{L_n(i,j)}(t)-a_{L_n(i,j)}(s)) \one_{\Delta_n}(i,j)\right)_{i,j=1}^N$ is the matrix-valued stochastic process with input sequence $ \{(a_m(t) - a_m(s))\}_{m \in \mathbb{Z}^d}$. 
	Now, using the trace formula of matrices, we get 
	\begin{align} \label{eqn:tracemultply}
		\Tr[ (A'_n(t,s) )^{p}] & = \sum_{\pi: \text{circuit} \atop \text{ of length } p} a'_{\pi}(t,s) \prod_{ i=1}^p \one_{\Delta_n}(\pi(i-1),\pi(i)), \\
		\Tr [(A'_n(t,s))^d (A_n(s))^{p-d}] & =  \sum_{\pi=\pi' \cup \pi'': \text{circuit} \atop \text{ of length } p} a'_{\pi'}(t,s) a_{\pi''}(s)\prod_{ i=1}^p \one_{\Delta_n}(\pi(i-1),\pi(i)),  \nonumber
	\end{align}
	where for a circuit $\pi:\{0,1,2,\ldots, p\} \rightarrow \{1,2,\ldots, N\}$,
	\begin{align*}
		a'_{\pi}(t,s) & = [a_{L(\pi(0),\pi(1))}(t)-a_{L(\pi(0),\pi(1))}(s)] 
		\cdots [a_{L(\pi(p-1),\pi(p))}(t) - a_{L(\pi(p-1),\pi(p))}(s)], \\
		a'_{\pi'}(t,s) & = [a_{L(\pi(0),\pi(1))}(t)-a_{L(\pi(0),\pi(1))}(s)] 
		\cdots [a_{L(\pi(d),\pi(d+1))}(t) - a_{L(\pi(d),\pi(d+1))}(s)], \\
		a_{\pi''}(s) & = a_{L(\pi(d+1),\pi(d+2))}(s) \cdots a_{L(\pi(p-1),\pi(p))}(s).
	\end{align*}
	Note that for a fixed $\pi$ and $i \in \{0,1, \ldots, p\}$, 
	$$a_{L(\pi(i),\pi(i+1))}(t)-a_{L(\pi(i),\pi(i+1))}(s) \distas{D} a_{L(\pi(i),\pi(i+1))}(t,s) \mbox{ and }  a_j(t) \distas{D} N(0, t).$$ 
	Therefore, $$x_{L(\pi(i),\pi(i+1))}  \distas{D} \frac{1}{\sqrt{
			t-s}} a_{L(\pi(i),\pi(i+1))}(t,s) \mbox{ and } y_{L(\pi(i),\pi(i+1))}  \distas{D} \frac{1}{\sqrt{s}} a_{L(\pi(i),\pi(i+1))}(s),$$
	where $\{x_{i}\}$ and $\{y_{i}\}$ are independent Gaussian random variables with mean zero and variance one. 
		Now on substituting (\ref{eqn:tracemultply}) with above notation in (\ref{eqn:trace t-s}), we get 
	\begin{align} \label{eqn:Z_p}
		\Tr (A_n(t))^{p}- \Tr (A_n(s))^{p}
		& \distas{D}   \sum_{\pi} \Big( (t-s)^{\frac{p}{2}} x_{\pi} + \sum_{d=1}^{p-1} \binom{p}{d}  (t-s)^{\frac{d}{2}} {s}^{\frac{p-d}{2}} x_{\pi'}  y_{\pi''}  \Big) \prod_{ i=1}^p \one_{\Delta_n}(\pi(i-1),\pi(i))  \nonumber \\
		& = \sqrt{t-s}  \sum_{\pi} Z_{\pi} I(\pi, \Delta_n,p), \mbox{ say},
	\end{align}  
    \vspace{-4mm}
	\begin{align*}
		& \text{where }Z_{\pi}=(t-s)^{\frac{p-1}{2}} x_{\pi} + \sum_{d=1}^{p-1} \binom{p}{d}  (t-s)^{\frac{d-1}{2}} {s}^{\frac{p-d}{2}} x_{\pi'}  y_{\pi''}, \\
         &I(\pi, \Delta_n,p)= \prod_{ i=1}^p \one_{\Delta_n}(\pi(i-1),\pi(i)),  \quad x_{\pi} = x_{L(\pi(0),\pi(1))} \cdots x_{L(\pi(p-1),\pi(p))}, \\&  x_{\pi'}  y_{\pi''} =  x_{L(\pi(0),\pi(1))} \cdots x_{L(\pi(d),\pi(d+1))}  y_{L(\pi(d+1),\pi(d+2))} \cdots y_{L(\pi(p-1),\pi(p))}.
	\end{align*}
	Finally, by substituting (\ref{eqn:Z_p}) in (\ref{eqn:w_p(t-s)}), we get
	\begin{align} \label{eqn:w_p^4}
		\E[\kap_{p}(t) -\kap_{p}(s)]^4 & = \frac{(t-s)^2}{N^{2p+2}} \sum_{\pi_1, \pi_2, \pi_3, \pi_4} \E \Big[ \prod_{k=1}^4 (Z_{\pi_k}- \E Z_{\pi_k})  \Big] \prod_{k=1}^4 I(\pi_k, \Delta_n,p).
	\end{align}
	Now to simplify \eqref{eqn:w_p^4}, we define the following notion of matchings of circuits.
    \begin{definition}
       For a  link function $L$ and a circuit $\pi$, define the multi-set $\operatorname{Im}_\pi(L):=\{L((\pi(i-1),\pi(i)): 1 \leq i \leq p\}$. For circuits $\pi_1, \pi_2, \ldots, \pi_k$, we say
        \begin{enumerate}
            \item [(i)] $\{\pi_1, \pi_2, \ldots, \pi_k\}$ are \textit{jointly-matched} if each element of the additive union $\uplus \operatorname{Im}_{\pi_i}(L)$ appears at least twice.
            \item [(ii)]$\{\pi_1, \pi_2, \ldots, \pi_k\}$ are said to be \textit{cross-matched} if for all $i$, there exists $j \neq i$ such that $\operatorname{Im}_{\pi_i}(L) \cap \operatorname{Im}_{\pi_j}(L) \neq \emptyset$.
        \end{enumerate}
    \end{definition}
	Consider $\{\pi_1, \pi_2, \pi_3, \pi_4 \}$, and suppose there exists an $i$ such that $\operatorname{Im}_{\pi_i}(L) \cap \operatorname{Im}_{\pi_j}(L) = \emptyset$ for all $j \neq i$, then from the independence of the entries, we have $ \E \Big[ \prod_{k=1}^4 (Z_{\pi_k}- \E Z_{\pi_k})  \Big]=0$. Thus, this case has zero contribution in  (\ref{eqn:w_p^4}). 
	Hence 
	\begin{align*} 
		\E[\kap_{p}(t) -\kap_{p}(s)]^4 & = \frac{(t-s)^2}{N^{2p+2}} \sum_{B_{p,4}} \E \Big[ \prod_{k=1}^4 (Z_{\pi_k}- \E Z_{\pi_k})  \Big] \prod_{k=1}^4 I(\pi_k, \Delta_n,p),
	\end{align*}
	where 
	$
		B_{p,4} = \{(\pi_1, \pi_2, \pi_3, \pi_4): \{\pi_1, \pi_2, \pi_3, \pi_4\} \mbox{ are jointly-matched and cross-matched}\}$.
	Since the moments of $\{x_i\},\{ y_i\}$ are uniformly bounded and $|I(\pi_k, \Delta_n,p)| \leq 1$, we have 
	\begin{align*} 
		\E[\kap_{p}(t) -\kap_{p}(s)]^4  \leq  \frac{(t-s)^2}{N^{2p+2}} C_T \times \# B_{p,4}, 
	\end{align*}
	where $C_T$ is a constant depending only on $p,T$ and independent of $n$.
	Now observe from Lemma 1.4.3 (a) of \cite{bose_patterned} that $\#B_{p,4} = O(N^{2p+2})$. Therefore from the above inequality, we get
	\begin{align*} 
		\E[\kap_{p}(t) -\kap_{p}(s)]^4 & \leq M_T (t-s)^2 \ \ \ \forall \ n \in \mathbb{N} \ \mbox{and } t,s \in[0,\ T],
	\end{align*}
	where $M_T$ is a positive constant depending only on $p, T$.
	This completes the proof of Proposition \ref{pro:tight} with $\alpha = 4$ and $\beta = 1$.
\end{proof}
Therefore by Lemma \ref{lem:two conditions} and Lemma \ref{res:exist of process}, we get for each $p$,
\begin{equation*}
	\big\{	\kap_{p}(t) ; t \geq 0\big\}  \stackrel{\mathcal{D}}{\rightarrow} \big\{ \Gamma_p(t) ; t \geq 0\big\},
\end{equation*}
where for each $r$, $(\Gamma_p(t_1),  \ldots, \Gamma_p(t_r))  \distas{D} \Gamma_{t_1,t_2,\ldots , t_r}$ with $\Gamma_{t_1,t_2,\ldots , t_r}$ given in Lemma \ref{lem:multi_dim}.

\subsection{Gaussian process} \label{subsec:gaussian}
In this subsection, we prove the last assertion of Theorem \ref{thm:process_brownian}, i.e., if the limiting distribution of $\kappa_p$ is Gaussian, say $\Gamma_{p}$,  then the limiting process  $ \big\{ \Gamma_{p}(t) ; t \geq 0\big\}$ is a Gaussian process.

\begin{proof}
	To show that $ \big\{ \Gamma_{p}(t) ; t \geq 0\big\}$ is a Gaussian process, it is enough to show that for any fixed $t_1, t_2, \ldots, t_r$;
 $(\Gamma_{p}(t_1), \Gamma_{p}(t_2), \ldots, \Gamma_{p}(t_r))$ is a 
 Gaussian vector. For this, we use  Wick's formula and the M-determinacy of Gaussian distribution. We prove that
 \begin{align}
     \lim_{n \to \infty}\mathbb{E}[	\kap_p^{k_1}\left(t_1\right)  \cdots \kap_p^{k_r}\left(t_r\right)]
     &= \lim_{n \to \infty} \sum_{\Pi \in \PP_2(2k)} \prod_{\{u,v\} \in \Pi}\E[\kappa_p(s_u)\kappa_p(s_v)] \label{eqn:Gaussian_wick},
 \end{align}
for all $k_1,k_2,\ldots, k_r$ with $k_1+\cdots +k_r=2k$. Here, $s_i=t_m$ whenever $k_1+k_2+\cdots +k_{m-1}<i \leq k_1+k_2+\cdots+ k_m$. Note that \eqref{eqn:Gaussian_wick} holds trivially when $k_1+k_2+\cdots + k_r$ is odd.

 Now, using \eqref{eqn: E eta_p time product} and condition (iii) of Theorem \ref{thm:process_brownian}, we obtain
 \begin{align*}
     \E [\Gamma_p^2] &= \lim_{n \to \infty} \frac{1}{N^{p+1}} \sum_{W \in SP(p,2)} \sum_{\Pi_{L_n,\Delta_n}(W)}  \mathbb{E}\left(\prod_{i=1}^2 \big(a_{\pi_i}-\mathbb{E} a_{\pi_i} \big)\right)\\
     &= \lim_{n \to \infty}  \sum_{W \in SP(p,2)} \frac{\# \Pi_{L_n, \Delta_n}(W)}{N^{p+1}}  \mathbb{E}\left(\prod_{i=1}^2 \big(a_{w_i}-\mathbb{E} a_{w_i} \big)\right).
 \end{align*}
This implies that
\begin{align}
    \bigl(\mathbb{E}[\Gamma_p^2]\bigr)^k
    &= \lim_{n\to\infty}
    \Bigg(
        \sum_{W\in SP(p,2)}
        \frac{\#\Pi_{L_n,\Delta_n}(W)}{N^{p+1}}
        \mathbb{E}\Bigg[
            \prod_{j=1}^{2}
            \bigl(a_{w_j}-\mathbb{E}[a_{w_j}]\bigr)
        \Bigg]
    \Bigg)^k \nonumber\\
    &= \lim_{n\to\infty}
    \sum_{W_1,\ldots,W_k\in SP(p,2)}
    \frac{\#\Pi_{L_n,\Delta_n}(\widetilde{W})}{N^{pk+k}}
    \prod_{i=1}^{k}
    \mathbb{E}\Bigg[
        \prod_{j=1}^{2}
        \bigl(a_{w_i^{(j)}}-\mathbb{E}[a_{w_i^{(j)}}]\bigr)
    \Bigg],
    \label{eqn:Gamma_second_moment}
\end{align}
where $W_i=(w_i^{(1)},w_i^{(2)})$ for $1\leq i\leq k$ and the letters of $W_i$ are represented by $a^{(i)},b^{(i)},\ldots$.
Here $\widetilde{W}=(w_1^{(1)},w_1^{(2)},\ldots,w_k^{(1)},w_k^{(2)})$ and the replacement of
$\prod_{i=1}^{k}\#\Pi_{L_n,\Delta_n}(W_i)$ by
$\#\Pi_{L_n,\Delta_n}(\widetilde{W})$ in the second line is justified by
observation (i) in the proof of Theorem~1.1 of
\cite{KMS_UnifiedLES_24}.
Note that the set of clusters of $\widetilde{W}$ is $\{\{1,2\},\{3,4\},\ldots , \{2k-1,2k\}\}$. Thus, the sum in \eqref{eqn:Gamma_second_moment} ranges only over
sentences in $SP(p,2k)$ with this prescribed cluster partition.

Let $\widetilde{SP}(p,2k)$ denote the set of all sentences
$W\in SP(p,2k)$ such that
$W_C \in SP(p,2)$ for all clusters $C$ of $W$.
For $W=(w_1,\ldots,w_{2k})\in\widetilde{SP}(p,2k)$, list its
clusters in the order of their smallest elements and write the $r$-th cluster as
$\{r_1,r_2\}$, with $r_1<r_2$.
Define $W'=(w_1',\ldots,w_{2k}')$ by
\[
    w_{2r-1}'=w_{r_1},
    \qquad
    w_{2r}'=w_{r_2},
    \qquad 1\leq r\leq k.
\]
Then $W'$ has the form considered in
\eqref{eqn:Gamma_second_moment}.
Further, each such sentence $W'$ has exactly $\#\PP_2(2k)$ preimages under
this map, one for each pair partition of $\{1,\ldots,2k\}$.
Since the value of the summand in \eqref{eqn:Gamma_second_moment} is the same for $W'$ and all its preimages, and since the product could be taken inside the expectation due to independence, we obtain
\begin{equation}
    \bigl(\mathbb{E}[\Gamma_p^2]\bigr)^k
    =
    \lim_{n\to\infty}
    \frac{1}{\#\PP_2(2k)}
    \sum_{W\in\widetilde{SP}(p,2k)}
    \frac{\#\Pi_{L_n,\Delta_n}(W)}{N^{pk+k}}
    \mathbb{E}\Bigg[
        \prod_{j=1}^{2k}
        \bigl(a_{w_j}-\mathbb{E}[a_{w_j}]\bigr)
    \Bigg].\label{eqn:Gamma_second_moment_final}
\end{equation}
Now by \cite[Proposition 3.2]{KMS_UnifiedLES_24},
\begin{align}
    \E[\Gamma_p^{2k}] &= \lim_{n \to \infty}  \sum_{W \in SP(p,2k)} \frac{\#\Pi_{L_n,\Delta_n}(W)}{N^{pk+k}}  \mathbb{E}\left(\prod_{i=1}^{2k} \big(a_{w_i}-\mathbb{E} a_{w_i}\big)\right) \nonumber\\
    &= \lim_{n \to \infty}  \sum_{W \in \widetilde{SP}(p,2k)} u(W)+  \lim_{n \to \infty} \sum_{W \in SP(p,2k) \setminus \widetilde{SP}(p,2k)} u(W),\label{eqn:Gamma_2k_moment}
\end{align}
where $u(W)= \frac{\#\Pi_{L_n,\Delta_n}(W)}{N^{pk+k}}  \mathbb{E}\left(\prod_{i=1}^{2k} \big(a_{w_i}-\mathbb{E} a_{w_i}\big)\right)$. Note that by Wick's formula, \eqref{eqn:Gamma_second_moment_final} and \eqref{eqn:Gamma_2k_moment} we have 
\begin{equation}
 \sum_{W \in SP(p,k) \setminus \widetilde{SP}(p,k)} u(W)=o(1).
 \label{eqn:residue_Wick}
\end{equation}
We now make the following claim.
\begin{claim}\label{claim:contri}
In this case, for  $W \in SP(p,2k) \setminus \widetilde{SP}(p,2k) $ such that all clusters of $W$ have even size, $\frac{\#\Pi_{L_n,\Delta_n}(W)}{N^{pk+k}}= o(1)$.
\end{claim}
\begin{proof}[Proof of claim.]
Consider a sentence $W \in SP(p,2k) \setminus \widetilde{SP}(p,2k)$ such that all clusters of $W$ are of even size. By \eqref{eqn:residue_Wick}, it follows that $u(W) \to 0$, as $n \to \infty$.
Further, by \eqref{eqn:contri_clusters} we have $\mathbb{E}\left(\prod_{i=1}^k \big(a_{w_i}-\mathbb{E} a_{w_i}\big)\right)$ is a non-zero constant. As a result, we conclude the existence of a constant $c>0$ such that 
\(
 u(W) \geq  c \times \frac{\#\Pi_{L_n,\Delta_n}(W)}{N^{pk+k}}.
\)
This implies that for sentences $W$ with all clusters of even size, $\frac{\#\Pi_{L_n,\Delta_n}(W)}{N^{pk+k}}=o(1)$, proving the claim.
\end{proof}
Now, we go back to the proof of the main assertion. Recall that 
\begin{equation*}
 	\mathbb{E}[	\kap_p^{k_1}\left(t_1\right)  \cdots \kap_p^{k_r}\left(t_r\right)] 
 	= \frac{1}{N^{pk+k}}  \sum_{W \in SP(p,2k)} \sum_{\Pi_{L_n,\Delta_n}(W)}  \mathbb{E}\left(\prod_{i=1}^{2k} \big(a_{w_i} (s_i)-\mathbb{E} a_{w_i}(s_i) \big)\right).
 \end{equation*}
If $W_C$ is a clique sentence of odd size $k$, then by \eqref{eqn:contri_clusters},
\[
\mathbb{E}\left(\prod_{i=1}^{2k} \big(a_{w_i} (s_i)-\mathbb{E} a_{w_i}(s_i) \big)\right)=0,
\] 
and in other cases $ \mathbb{E}\left(\prod_{i=1}^{2k} \big(a_{w_i} (s_i)-\mathbb{E} a_{w_i}(s_i) \big)\right)$ is bounded. Therefore by Claim \ref{claim:contri} it follows that for any choice of $0 < t_1 <t_2 <\cdots < t_r$, 
\begin{align*}
\lim_{n \to \infty}\mathbb{E}[	&\kap_p^{k_1}\left(t_1\right)  \cdots \kap_p^{k_r}\left(t_r\right)] = \lim_{n \to \infty}\frac{1}{N^{pk+k}}  \sum_{W \in \widetilde{SP}(p,2k)} \sum_{\Pi_{L_n,\Delta_n}(W)}  \mathbb{E}\left(\prod_{i=1}^{2k} \big(a_{w_i} (t_i)-\mathbb{E} a_{w_i}(t_i) \big)\right). 
\end{align*}
We rewrite the summation over $W \in \widetilde{SP}(p,2k)$ into summations over sentences $W_1,\ldots, W_k \in SP(p,2)$ (like in \eqref{eqn:Gamma_second_moment}) noting that each combination of $(W_i; 1 \leq i \leq k)$ corresponds to $\# \PP_2(2k)$ words in $\widetilde{SP}(p,2k)$. But here, unlike \eqref{eqn:Gamma_second_moment_final}, the contribution of each pair-partition is different, and we get that the above equation is equal to
\begin{align*}
 \lim_{n\to\infty}
    &\sum_{\Pi \in \PP_2(2k)}\sum_{W_1,\ldots,W_k\in SP(p,2)}  
    \frac{\prod_{i=1}^k\#\Pi_{L_n,\Delta_n}(W_i)}{N^{pk+k}}  \prod_{i=1}^{k}
    \mathbb{E}\Bigg[
        \prod_{j=1}^{2}
        \bigl(a_{w_i^{(j)}}(s_{i_1})-\mathbb{E}[a_{w_i^{(j)}}(s_{i_2})]\bigr)
    \Bigg],
\end{align*}
where $W_i=(w_i^{(1)},w_i^{(2)})$, $\{i_1,i_2\}$ is the $i$-th block of $\Pi$ and $\{s_i\}$ is as defined in \eqref{eqn:si}.
Now note that the above equation is exactly the right hand-side of \eqref{eqn:Gaussian_wick}, and this completes the proof.
 \end{proof}

\section{Proof of Theorem \ref{thm:process_brownian}: Ornstein-Uhlenbeck process entries} \label{sec:thm_OH}
In this section we prove Theorem \ref{thm:process_brownian} for Ornstein-Uhlenbeck process entries. 
The broad idea of the proof is similar to the proof for Brownian motion entries, and so we highlight the key differences.
\vskip2pt
\noindent\textbf{Finite-dimensional convergence.} 
Consider $\kap_p(t)$ from \eqref{eqn:kappa_p(t)} with the input  of $A_n(t)$ as independent OU processes. Note that in this case as well, \eqref{eqn: E eta_p time product} holds 
and further steps in the proofs in Section \ref{subsec:finite_dim} also go through. As a result, we get that for $W \notin SP(p,k)$, their contribution asymptotically goes to zero. For $W \in SP(p,k)$ and $\{s_i \}$ as defined by \eqref{eqn:si}, a direct calculation gives
\begin{align*}
	 	&\ \mathbb{E}\left[\prod_{i  \in C} (a_{w_i} (s_i)-\mathbb{E} a_{w_i}(s_i))\right]\\
        &=
	 \begin{cases}
	   \left(\frac{\sigma^2}{2\theta}\right)^{\frac{p\ell}{2}}\sum_{\pi \in \PP_2(\ell)} \prod_{\{r,b\} \in \pi}e^{-\theta|s_{i_r}-s_{i_b}|}
       & \text{ if } W_C=(w_{i_1},\ldots, w_{i_\ell}) \text{ is a clique sentence},\\
	 		\left(\frac{\sigma^2}{2\theta}\right)^{p} e^{-c_0\theta|s_{i_1}-s_{i_2}|} &	\text{ if } W_C=(w_{i_1},w_{i_2}) \in \PP_{2}(p,p), \\
	 		\left(\frac{\sigma^2}{2\theta}\right)^{p}(1+2e^{-2\theta|s_{i_1}-s_{i_2}|})
            &	\text{ if } W_C=(w_{i_1},w_{i_2}) \in \PP_{2,4}(p,p),
	 \end{cases}	 
\end{align*}
where for $(w_{i_1},w_{i_2}) \in \PP_{2}(p,p)$, $c_0$ is the number of  elements in $S_{w_{i_1}} \cap S_{w_{i_2}}$. This implies that $\mathbb{E}[	\kap_p^{k_1}\left(t_1\right) \kap_p^{k_2}\left(t_2\right) \cdots \kap_p^{k_r}\left(t_r\right)]$ converges for all non-negative choices of $k_1,k_2,\ldots,k_r$. The rest of the arguments in Section \ref{subsec:finite_dim} follow as Lemma \ref{lem:multi_dim} is applicable in this case as well.

\vskip2pt
\noindent
\textbf{Tightness.}  The idea and most of the arguments for establishing tightness for OU process are similar to the proof of Proposition \ref{pro:tight}. Here we outline only the main steps.

Since $\{a(t)\}$ is an OU process,
$\var(a(t))=\frac{\sigma^2}{2\theta}$ for all $t\geq 0  \mbox{ and }
\Cov(a(t),a(s))
=\frac{\sigma^2}{2\theta}
e^{-\theta|t-s|}$, for all $t,s\geq 0$.
Consequently,
\begin{equation*}
 \mathbb{E}[a(t)-a(s)]
=0  \quad \mbox{and} \quad 
\operatorname{Var}\big(a(t)-a(s)\big)
= \frac{\sigma^2}{\theta} \left(1-e^{-\theta|t-s|}\right).   
\end{equation*}
Noting that both $a(t)$ and $a(s)$ have normal distribution, it follows that
$a(t)-a(s) \overset{D}{\sim}\mathcal{N}\left(0,\, \frac{\sigma^2}{\theta}
\left(1-e^{-\theta|t-s|}\right)\right).$
Equivalently, if $Z\overset{D}{\sim} \NN(0,1)$, then
\begin{align}
    a(t)-a(s)
&\overset{D}{\sim} \sqrt{\frac{\sigma^2}{\theta} \left(1-e^{-\theta|t-s|}\right)}\,Z = \sqrt{c_{t,s}} Z, \mbox{ say},
\end{align}
where 
\( c_{t,s}=\frac{\sigma^2}{\theta}\left(1-e^{-\theta|t-s|}\right)
\). 
As a result, the equivalent of (\ref{eqn:Z_p}) in this case is
\begin{align} \label{eqn:Z_p_ou}
		\Tr (A_n(t))^{p}- \Tr (A_n(s))^{p}
		& \distas{D}   \sum_{\pi} \Big( (c_{t,s})^{\frac{p}{2}} x_{\pi} + \sum_{d=1}^{p-1} \binom{p}{d}  (c_{t,s})^{\frac{d}{2}} { \Big(\frac{\sigma^2}{2\theta}\Big)}^{\frac{p-d}{2}} x_{\pi'}  y_{\pi''}  \Big) \prod_{ i=1}^p \one_{\Delta_n}(\pi(i-1),\pi(i))  \nonumber \\
		& = \sqrt{c_{t,s}}  \sum_{\pi} Z_{\pi} I(\pi, \Delta_n,p), \mbox{ say},
	\end{align} 
	where $Z_\pi, x_{\pi}, x_{\pi'}$ and $y_{\pi''}$ are as in \eqref{eqn:Z_p}.
    Note that \eqref{eqn:Z_p_ou} is the same as \eqref{eqn:Z_p}, with the only differences being $c_{t,s}$ in place of $(t-s)$, and $x_{\pi'}$ and $y_{\pi''}$ are dependent in this case. Thus
    \begin{align} \label{eqn:w_p^4_OU}
		\E[\kap_{p}(t) -\kap_{p}(s)]^4 & = \frac{(c_{t,s})^2}{N^{2p+2}} \sum_{\pi_1, \pi_2, \pi_3, \pi_4} \E \Big[ \prod_{k=1}^4 (Z_{\pi_k}- \E Z_{\pi_k})  \Big] \prod_{k=1}^4 I(\pi_k, \Delta_n,p).
	\end{align}    
Note that, for $x>0$, $1-e^{-x}\leq x$ and therefore we obtain $1-e^{-\theta |t-s|} \leq \theta |t-s|$, and thus
\begin{align*}
\sqrt{c_{t,s}}
\leq\sqrt{ \frac{\sigma^2}{\theta} \left(\theta |t-s|\right)}
= \sigma\sqrt{|t-s|}.
\end{align*} 
Using the fact that the dependence between $x_{\pi'}$ and $y_{\pi''}$ does not affect the calculations based on matching, and noting that all fixed joint moments of $\{x_i\},\{ y_i\}$ are uniformly bounded,  we have 
	\begin{align*} 
		\E[\kap_{p}(t) -\kap_{p}(s)]^4 \leq \sigma^4 \frac{(t-s)^2}{N^{2p+2}} C_T' \times \# B_{p,4} \leq M_T' (t-s)^2 \ \ \ \forall \ n \in \mathbb{N} \ \mbox{and } t,s \in[0,\ T],
	\end{align*}
	where $M_T'$ is a positive constant depending only on $p, T$.
	This completes the proof of tightness part for OU process. \\ 
\noindent \textbf{Acknowledgment:} The research work of S.N. Maurya is supported by the
Inspire Faculty Fellowship (DST/INSPIRE/04/2024/003610). 
\bibliography{shambhubib}
\bibliographystyle{amsplain}
\end{document}